\documentclass[12pt]{amsart}
\usepackage{amsmath}
\usepackage{amsfonts}
\begin{document}
\numberwithin{equation}{section}

\def\1#1{\overline{#1}}
\def\2#1{\widetilde{#1}}
\def\3#1{\widehat{#1}}
\def\4#1{\mathbb{#1}}
\def\5#1{\frak{#1}}
\def\6#1{{\mathcal{#1}}}

\newcommand{\de}{\partial}
\newcommand{\R}{\mathbb R}
\newcommand{\al}{\alpha}
\newcommand{\tr}{\widetilde{\rho}}
\newcommand{\tz}{\widetilde{\zeta}}
\newcommand{\tv}{\widetilde{\varphi}}
\newcommand{\hv}{\hat{\varphi}}
\newcommand{\tu}{\tilde{u}}
\newcommand{\tF}{\tilde{F}}
\newcommand{\debar}{\overline{\de}}
\newcommand{\Z}{\mathbb Z}
\newcommand{\C}{\mathbb C}
\newcommand{\Po}{\mathbb P}
\newcommand{\zbar}{\overline{z}}
\newcommand{\G}{\mathcal{G}}
\newcommand{\So}{\mathcal{S}}
\newcommand{\Ko}{\mathcal{K}}
\newcommand{\U}{\mathcal{U}}
\newcommand{\B}{\mathbb B}
\newcommand{\oB}{\overline{\mathbb B}}
\newcommand{\Cur}{\mathcal D}
\newcommand{\Dis}{\mathcal Dis}
\newcommand{\Levi}{\mathcal L}
\newcommand{\SP}{\mathcal SP}
\newcommand{\Sp}{\mathcal Q}
\newcommand{\Ma}{\mathcal M}
\newcommand{\Co}{\mathcal C}
\newcommand{\Hol}{{\sf Hol}(\mathbb D, D)}
\newcommand{\Aut}{{\sf Aut}(\mathbb D)}
\newcommand{\D}{\mathbb D}
\newcommand{\oD}{\overline{\mathbb D}}
\newcommand{\oX}{\overline{X}}
\newcommand{\loc}{L^1_{\rm{loc}}}
\newcommand{\la}{\langle}
\newcommand{\ra}{\rangle}
\newcommand{\thh}{\tilde{h}}
\newcommand{\N}{\mathbb N}
\newcommand{\kd}{\kappa_D}
\newcommand{\Hr}{\mathbb H}
\newcommand{\ps}{{\sf Psh}}

\newcommand{\subh}{{\sf subh}}
\newcommand{\harm}{{\sf harm}}
\newcommand{\ph}{{\sf Ph}}
\newcommand{\tl}{\tilde{\lambda}}
\newcommand{\ts}{\tilde{\sigma}}

\def\v{\varphi}
\def\Re{{\sf Re}\,}
\def\Im{{\sf Im}\,}

\def\dist{{\rm dist}}
\def\const{{\rm const}}
\def\rk{{\rm rank\,}}
\def\id{{\sf id}}
\def\aut{{\sf aut}}
\def\Aut{{\sf Aut}}
\def\CR{{\rm CR}}
\def\GL{{\sf GL}}
\def\U{{\sf U}}

\def\la{\langle}
\def\ra{\rangle}

\newtheorem{theorem}{Theorem}[section]
\newtheorem{lemma}[theorem]{Lemma}
\newtheorem{proposition}[theorem]{Proposition}
\newtheorem{corollary}[theorem]{Corollary}

\theoremstyle{definition}
\newtheorem{definition}[theorem]{Definition}
\newtheorem{example}[theorem]{Example}

\theoremstyle{remark}
\newtheorem{remark}[theorem]{Remark}
\numberwithin{equation}{section}

\title[Aleksandrov-Clark measures and semigroups]{Aleksandrov-Clark measures and semigroups of analytic functions in the unit disc}
\author[F. Bracci]{Filippo Bracci$^\dag$}
\address{F. Bracci: Dipartimento Di Matematica\\
Universit\`{a} di Roma \textquotedblleft Tor Vergata\textquotedblright\ \\
Via Della Ricerca Scientifica 1, 00133 \\
Roma, Italy} \email{fbracci@mat.uniroma2.it}
\thanks{}
\date{\today}
\author[M.D. Contreras]{Manuel D. Contreras$^\ddag$}
\address{M.D. Contreras \and  S. D\'{\i}az-Madrigal: Camino de los Descubrimientos, s/n\\
Departamento de Matem\'{a}tica Aplicada II \\
Escuela T\'{e}cnica Superior de Ingenieros\\
Universidad de Sevilla\\
41092, Sevilla\\
Spain.} \email{contreras@us.es, madrigal@us.es}
\author[S. D\'{\i}az-Madrigal]{Santiago D\'{\i}az-Madrigal$^\ddag$}

\thanks{$^\dag$Partially supported by {\sl Progetto PRIN
Azioni di gruppi su variet\`a CR e complesse, spazi di moduli,
teoria geometrica delle funzioni e dinamica olomorfa}.\\
\indent $^\ddag$Partially supported by the \textit{Ministerio
de Ciencia y Tecnolog\'{\i}a} and the European Union (FEDER)
projects BFM2003-07294-C02-02 and MTM2006-14449-C02-01 and by \textit{La Consejer\'{\i}a
de Educaci\'{o}n y Ciencia de la Junta de Andaluc\'{\i}a.}}

\begin{abstract}
In this paper we prove a formula describing the infinitesimal
generator of a continuous semigroup $(\v_t)$ of holomorphic
self-maps of the unit disc with respect to a boundary regular
fixed point. The result is based on Alexandrov-Clark measures
techniques. In particular we prove that the Alexandrov-Clark
measure  of $(\v_t)$ at a boundary regular fixed points is
differentiable (in the weak$^\ast$-topology) with respect to
$t$.
\end{abstract}

\maketitle

\section{Introduction}

The aim of the present note is to study the incremental ratio
of Aleksandrov-Clark measures (sometimes called spectral
measures) of continuous semigroups of the unit disc at boundary
regular fixed points, obtaining a measure-theoretic
generalization of the well renowned Berkson-Porta formula at
the Denjoy-Wolff point.

To state our results, we briefly recall the notion of
Aleksandrov-Clark measures and semigroups as needed for our
aims (for details on Aleksandrov-Clark measures we refer the
reader to the recent surveys \cite{Mat-Ste}, \cite{Pol-Sar},
\cite{Sak} and the references therein; while we refer to
\cite{Abate} and \cite{Shoiket} for more about iteration theory
and  semigroups).

Let $\D:=\{z\in \C: |z|<1\}$ be the unit disc. Let $f:\D\to \D$
be holomorphic. Fix $\tau\in \de\D$ and consider the positive
harmonic function $\Re (\tau+f(z))/(\tau-f(z))$. Then there exists a
non-negative and finite Borel measure $\mu_{f,\tau}$ (called
the {\sl Aleksandrov-Clark measure} of $f$ at $\tau$) on $\de
\D$ such that
\begin{equation}\label{AleCl}
\Re \frac{\tau+f(z)}{\tau-f(z)}=\int_{\de \D} P(\zeta,z)
d\mu_{f,\tau}(\zeta) \quad \hbox{for all\ } z\in \D,
\end{equation}
where $P(\zeta,z)=\frac{1-|z|^2}{|z-\zeta|^2}$ is the Poisson
kernel.

We recall that a point $\zeta \in \de \D$ is said to be a
boundary contact point for $f:\D \to \D$ if $\lim_{r\rightarrow 1}f(r
\zeta )=\tau \in \de \D.$ In such a case, as customary, we
write $f^\ast(\zeta):=\tau.$ It is a remarkable fact that  the
{\sl angular derivative} at a boundary contact point $\zeta$
always exists (possibly infinity). Namely, the following
non-tangential (or angular) limit exists  in the Riemann
sphere:
\[
f^{\prime}(\zeta):=\angle\lim\limits_{z\rightarrow\zeta}\frac{f(z)-f^\ast(\zeta)}{z-\zeta}.
\]
The modulus of $f^{\prime}(\zeta)$ is known as the
boundary dilatation coefficient of $f$ at $\tau$.
A point $\tau\in \de \D$ is  a {\sl boundary regular fixed
point}, BRFP for short, for $f$  if $f^\ast(\tau)=\tau$ and
$f'(\tau)$ is finite. If $\tau$ is a BRFP for $f$ then the
classical Julia-Wolff-Carath\'eodory theorem (see, {\sl e.g.},
\cite[Prop. 1.2.8, Thm. 1.2.7]{Abate}) asserts that the
non-tangential limit $\angle\lim_{z\to \tau} f'(z)=f'(\tau)$
and $f'(\tau)\in (0,+\infty)$.

In 1929, R. Nevanlinna obtained a very deep relationship
between   angular derivatives and Aleksandrov-Clark measures.
Namely,  (see, {\sl e.g.}, \cite[p. 61]{Abate}, \cite[Thm.
3.1]{Sak}) he proved

\begin{theorem}[Nevanlinna]\label{nevanlinna}
Let $f:\D\to\D$ be a holomorphic self-map and $\zeta \in\de\D$.
Then $\zeta$  is a boundary contact point of $f$ with
$f'(\zeta)\in\C$ if and only if for some $\tau \in \de \D$ the
Aleksandrov-Clark measure $\mu_{f,\tau}$ has an atom at $\zeta$
(that is, $\mu_{f,\tau}(\{\zeta\})>0$). In this case, it
follows $f^\ast(\zeta)=\tau$ and
$\mu_{f,\tau}(\{\zeta\})=1/|f'(\zeta)|$.
 \end{theorem}

A (continuous) semigroup $(\v_t)$ of holomorphic self-maps is a
continuous homomorphism from the additive semigroup of
non-negative real numbers and the composition semigroup of all
holomorphic self-maps of $\D$ endowed with the compact-open
topology. It is well known after the basic work of Berkson and
Porta \cite{Berkson-Porta} that in fact the dependence of the
semigroup $(\v_t)$ on the parameter $t$ is real-analytic and there
exists a holomorphic vector field $G:\D\to \C$, called the {\sl
infinitesimal generator of the semigroup}, such that
\begin{equation}\label{genera}
    \frac{\de \v_t}{\de t}(z)=G(\v_t(z))
\end{equation}
 for all $z\in \D$.

A point $\tau\in \de \D$ is a {\sl boundary regular fixed
point}, BRFP for short, for the semigroup $(\v_t)$ provided it
is a BRFP for $\v_t$ for all $t\geq 0$ (and this is the case if
and only if $\tau$ is a BRFP for $\v_t$ for some $t>0$, see
\cite{CoDiPo-Fennicae}). According to \cite[Thm. 1]{CoDiPo},
the point $\tau\in\de \D$ is a BRFP for $(\v_t)$ if and only if
$G(\tau)=0$ as non-tangential limit and the non-tangential
limit $\angle\lim_{z\to \tau} G'(z)=\lambda$ exists finitely.
Moreover, if $\tau$ is a BRFP for $(\v_t)$ then
$\v_t'(\tau)=e^{\lambda t}$
and $\lambda \in \R$. It is well known that if $(\v_t)$ has no
fixed point in $\D$ then there exists a unique boundary
regular fixed point $\tau\in\de\D$, called the {\sl
Denjoy-Wolff point} of the semigroup, such that $\v_t(z)\to
\tau$ as $t\to\infty$ for all $z\in\D$.

In the rest of the paper we will denote by $dm$ the Lebesgue
measure on $\de \D$ normalized so that $m(\de\D)=1$ and by
$\delta_\xi$ the Dirac atomic measure concentrated at
$\xi\in\de\D$.

Let $(\v_t)$ be a semigroup  and  $\tau\in\de\D$. We will
denote by $\mu_{t,\tau}$ the Aleksandrov-Clark measure of
$\v_t$ at $\tau$. It can be checked (test with Poisson kernels
and use the density of their span in $C(\de \D)$) that
$\{\mu_{t,\tau}\}$ is continuous (in the weak$^\ast$-topology)
with respect to $t$. We will prove that it is actually
differentiable at $0$. Indeed, our first result is the
following (note that $\mu _{0,\tau}=\delta _\tau$):

\begin{proposition}\label{mainproposition}
Let $(\v_t)$ be a continuous semigroup of holomorphic self-maps
of the unit disc $\D$. Let $\tau\in\de\D$ be a boundary regular
fixed point for $(\v_t)$ with boundary dilatation coefficients
$(e^{\lambda t})$. Then there exists a positive measure $\mu $
on $\de\D$ such that
\begin{equation}\label{maineq-intro}
\frac{\mu_{t,\tau}-\delta_\tau}{t}\stackrel{\hbox{w}^\ast}{\longrightarrow}
-\lambda \delta_\tau + \mu, \
 \quad \hbox{as } t\to 0.
\end{equation}
\end{proposition}

The measure $\mu$ in \eqref{maineq-intro} is strictly related
to the infinitesimal generator of the semigroup, see
Proposition \ref{Glimit}. From such a formula we will obtain
our main result:

\begin{theorem}\label{main}
Let $(\v_t)$ be a continuous semigroup of holomorphic self-maps
of the unit disc $\D$ and  let $G$ be its infinitesimal
generator. Let $\tau\in\de\D$ be a boundary regular fixed point
for $(\v_t)$  with boundary dilatation coefficients
$(e^{\lambda t})$. Then there exists a unique $p:\D\to\C$
holomorphic, with $\Re p\geq 0$ and $\angle \lim _{z\to
\tau}(z-\tau)p(z)=0$  such that
\begin{equation}\label{descom-geninf}
G(z)=({\bar \tau}z-1)(z-\tau)\left[p(z)-\frac{\lambda}{2}\frac{\tau+z}{\tau-z}\right]
\quad\hbox{for all}\ z\in\D.
\end{equation}
Conversely, given $p:\D\to\C$ holomorphic with $\Re p\geq 0$
and $\angle \lim _{z\to \tau}(z-\tau)p(z)=0$, $\tau \in \de \D$
and $\lambda \in \R$,
 the function $G$ defined as in \eqref{descom-geninf}
is the infinitesimal generator of a semigroup of holomorphic
self-maps of the unit disc for which $\tau$ is a boundary
regular fixed point with boundary dilatation coefficients
$(e^{\lambda t})$.
\end{theorem}

In particular, $\tau$ is the Denjoy-Wolff point of $(\v_t)$ if
and only if $\lambda \leq 0$ and, if this is the case,
\[
\Re\left(\frac{G(z)}{({\bar \tau}z-1)(z-\tau)}\right)\geq 0 \quad\hbox{for all}\ z\in\D,
\]
recovering in this way the celebrated Berkson-Porta
representation formula when $\tau$ belongs to $\de \D$ \cite{Berkson-Porta}.

The plan of the paper is the following. In the second section
 we compute the singular part of Aleksandrov-Clark measures
of $N$-to-$1$ mappings (in particular for $N=1$, univalent
maps). In the third section we    will use such computation to
prove Proposition \ref{mainproposition} and Theorem \ref{main}.
In the final section we discuss some consequences of our
results.

\section{Singular parts of Aleksandrov-Clark measures for $N$-to-$1$ mappings}

P.J. Nieminen and E. Saksman \cite[p. 3186]{Nie-Sak} already
remarked that for holomorphic $N$-to-$1$ self-maps of the unit
disc the singular part of the corresponding Aleksandrov-Clark
measures is discrete. In this section we enhance this result
explicitly computing such a singular part.

Given a positive Borel measure $\varrho$ on $\de \D$ we will
write $\varrho=\varrho^s+\varrho^a$  for its Lebesgue
decomposition in the singular part $\varrho^s$ and the absolutely
continuous part $\varrho^a$ with respect to the Lebesgue
measure.

\begin{proposition}\label{Na1}
Let $f:\D\to\D$ be a $N$-to-$1$ ($N\geq 1$) holomorphic   map
and let $\tau\in\de\D$. Then  there exist $0\leq m\leq N$ and
$\zeta_1,\ldots, \zeta_m\in\de \D$ such that
$f^\ast(\zeta_j)=\tau$, the non-tangential limit $f'(\zeta_j)$
of $f'$ at $\zeta_j$ exists finitely for $j=1,\ldots, m$ and
\begin{equation}\label{singolareNa1}
    \mu_{f,\tau}^s=\sum_{k=1}^m
    \frac{1}{|f'(\zeta_k)|}\delta_{\zeta_k}.
\end{equation}
Moreover, if $x\in\de \D\setminus\{\zeta_1,\ldots, \zeta_m\}$
is such that $f^\ast(x)=\tau$ then $\limsup_{z\to
x}|f'(z)|=\infty$.
\end{proposition}

In order to prove Proposition \ref{Na1} we need the following
lemma:

\begin{lemma}\label{lem1}
Let $f:\D\to\D$ be holomorphic and let $x,y\in\de\D$ be such
that $f^\ast(x)=y$. If $\{z_n\}\subset\D$ is a sequence
converging {\sl tangentially} to $x$ such that $\{f(z_n)\}$
converges {\sl non-tangentially} to $y$ then
\[
\lim_{n\to\infty}\frac{1-|f(z_n)|}{1-|z_n|}=+\infty.
\]
\end{lemma}
\begin{proof}
Since $\{f(z_n)\}$ converges non-tangentially to $y$ then there
exists $C>0$ such that for all $n\in \N$ it holds
\[
\frac{1-|f(z_n)|}{|y-f(z_n)|}\geq C.
\]
Therefore
\begin{equation*}
\begin{split}
\lim_{n\to\infty}\frac{1-|f(z_n)|}{1-|z_n|}&=\lim_{n\to \infty}
\frac{1-|f(z_n)|}{|y-f(z_n)|}\cdot
\frac{|y-f(z_n)|}{|x-z_n|}\cdot \frac{|x-z_n|}{1-|z_n|}\\ &\geq
C \cdot
\left[\liminf_{n\to\infty}\frac{|y-f(z_n)|}{|x-z_n|}\right]\cdot
\left[\lim_{n\to\infty}\frac{|x-z_n|}{1-|z_n|}\right]=+\infty,
\end{split}
\end{equation*}
since $\liminf_{n\to\infty}\frac{|y-f(z_n)|}{|x-z_n|}>0$ by
Julia's Lemma (see, {\sl e.g.}, \cite{Abate}) and
$\lim_{n\to\infty}\frac{|x-z_n|}{1-|z_n|}=+\infty$ by
hypothesis.
\end{proof}

\begin{proof}[Proof of Proposition \ref{Na1}]
According to \cite[Prop. 5.3]{Nie-Sak} there exists a sequence
$\{z_n\}\subset \D$ which converges non-tangentially to $\tau$
such that
\begin{equation}\label{Nisanumb}
    \|\mu_{f,\tau}^s\|=\lim_{n\to \infty} \frac{\sum_{w\in
    f^{-1}(z_n)}  \log|w|}{\log|z_n|}.
\end{equation}
Since $f$ is $N$-to-$1$, the number of preimages $f^{-1}(z_n)$
for each $n\in\N$ is (at most) $N$. Let
$M:=\limsup_{n\to\infty}\sharp\{f^{-1}(z_n)\}$ be the supremum
limit of the number of preimages of $z_n$.  Notice that $0\leq
M\leq N$. If $M=0$, namely, if $f^{-1}(z_n)=\emptyset$
eventually, then $\mu_{f,\tau}^s=0$. Assume that $M>0$. Then
there exists a subsequence $\{z_{n_k}\}$ of $\{z_n\}$ such that
$\sharp\{f^{-1}(z_{n_m})\}=M$. By \eqref{Nisanumb} the mass
$\|\mu_{f,\tau}^s\|$ is expressed as a limit, thus we can
replace $\{z_n\}$ with $\{z_{n_k}\}$ and assume directly that
$\sharp\{f^{-1}(z_{n})\}=M$ for all $n$. Let us denote by
$\{w_{n,1},\ldots, w_{n,M}\}$ the preimages of $z_n$. Again by
\eqref{Nisanumb}, up to extracting subsequences, we can assume
that $\{w_{n,j}\}$ is converging to $\sigma_j\in\oD$ for
$j=1,\ldots, M$.

Fix $j\in\{1,\ldots,M\}$. Clearly, since $f$ is open,
$\sigma_j\in\de\D$. Moreover, since
$\lim_{n\to\infty}\frac{\log|w_{n,j}|}{1-| w_{n,j} |}=1$ (and
similarly for $z_n$), it follows that
\[
\lim_{n\to\infty}\frac{\log|w_{n,j}|}{\log|z_{n}|}=\lim_{n\to\infty}\frac{1-|w_{n,j}|}{1-|z_{n}|}
=\lim_{n\to\infty}\frac{1-|w_{n,j}|}{1-|f(w_{n,j})|}.
\]
By Lemma \ref{lem1}, if $\{w_{n,j}\}$ converges to $\sigma_j$
tangentially then
$\lim_{n\to\infty}\frac{1-|w_{n,j}|}{1-|f(w_{n,j})|}=0$.
Suppose then that  $\{w_{n,j}\}$ converges to $\sigma_j$
non-tangentially. By the classical Julia-Wolff-Carath\'eodory
theorem, either
$\lim_{n\to\infty}\frac{1-|w_{n,j}|}{1-|f(w_{n,j})|}=0$ or the
non-tangential limit $f'(\sigma_j)$ of $f'$ at $\sigma_j$
exists finitely and
$\lim_{n\to\infty}\frac{1-|w_{n,j}|}{1-|f(w_{n,j})|}=|f'(\sigma_j)|^{-1}$.

Let $0\leq m\leq M$ and $\{\zeta_1,\ldots,
\zeta_m\}\subseteq\{\sigma_1,\ldots, \sigma_M\}$ be the biggest
possible subset such that the non-tangential limit
$f'(\zeta_j)$ exists finitely for every $j$. Then we rewrite (\ref{Nisanumb}) as
\[
\|\mu_{f,\tau}^s\|=\sum_{j=1}^{m}\frac{1}{|f'(\zeta_j)|}.
\]
Moreover, by Theorem \ref{nevanlinna}, we know that
\[
\sum_{j=1}^{m}\frac{1}{|f'(\zeta_j)|}\delta_{\zeta_j}\leq \mu_{f,\tau}^s,
\]
Then \eqref{singolareNa1} holds.

To conclude the proof, we need to show that if $x\in\de
\D\setminus\{\zeta_1,\ldots, \zeta_m\}$ is such that
$f^\ast(x)=\tau$ then $|f'(x)|=\infty$. Indeed,
if this were not true then by the Julia-Wolff-Carath\'eodory
theorem the non-tangential derivative $f'(x)$ would exist and
$|f'(x)|^{-1}\delta_x$ would be part of $\mu_{f,\tau}^s$ (see,
{\sl e.g.}, \cite[Thm. 3.1]{Sak}), contradicting
\eqref{singolareNa1}.
\end{proof}

\begin{corollary}\label{univaleNev}
Let $f:\D\to\D$ be a univalent map and let $\tau\in\de\D$. Then
either $\mu_{f,\tau}^s=0$ or there exists a unique point
$\zeta\in\de \D$ such that $f^\ast(\zeta)=\tau$, the
non-tangential limit $f'(\zeta)$ of $f'$ at $\zeta$ exists
finitely  and
\begin{equation}\label{singolareNa2}
    \mu_{f,\tau}^s=
    \frac{1}{|f'(\zeta)|}\delta_{\zeta}.
\end{equation}
Moreover, if $x\in\de \D\setminus\{\zeta\}$ is such that
$f^\ast(x)=\tau$ then $\limsup_{z\to x}|f'(z)|=\infty$.
\end{corollary}

\begin{remark}
Corollary \ref{univaleNev} implies in particular that for a
univalent self-map $f$ of the unit disc and any point
$\tau\in\de\D$ there exists at most one point $x\in\de\D$ such
that $f^\ast(x)=\tau$ and the non-tangential limit of $f'$
exists finitely at $x$. This latter fact can also be  proved
directly,  see \cite[Lemma 8.2]{CoPo}.
\end{remark}

\begin{corollary}\label{univalegrupNev}
Let $(\v_t)$ be a continuous semigroup of holomorphic self-maps
of $\D$. Suppose that $\tau\in\de\D$ is a BRFP for $(\v_t)$
with boundary dilatation coefficients $(e^{\lambda t})$. Then
\begin{equation}\label{singolareNa3}
    \mu_{t,\tau}^s=e^{-\lambda t}\delta_{\tau}.
\end{equation}
\end{corollary}

\begin{proof}
For every $t\geq 0$ the map $\v_t$ is univalent (see, {\sl
e.g.}, \cite{Abate} or \cite{Shoiket}). Therefore by
Corollary~\ref{univaleNev} it follows that $\mu_{t,\tau}^s=
    \frac{1}{|\v_t'(\tau)|}\delta_{\tau}$ and since
    $\v_t'(\tau)=e^{\lambda t}$, we are done.
\end{proof}

\begin{remark}
In the proof of Proposition \ref{Na1} and as a byproduct of Theorem \ref{nevanlinna},
we used that for an arbitrary holomorphic self-map of the unit disk $f$,
given $\tau \in \de \D$ and $\zeta _1, ... \zeta _n $ different points in $\de\D$ such that
$f^\ast (\zeta_j)=\tau$ with $f'(\zeta_j)\in\C$  for all $j=1,...,n$, then
\[
\sum_{j=1}^{n}\frac{1}{|f'(\zeta_j)|}\delta_{\zeta_j}\leq \mu_{f,\tau}^s.
\]
In particular, we have
\begin{equation}\label{Cowen-Pomm-Ine}
\sum_{j=1}^{n}\frac{1}{|f'(\zeta_j)|}\leq \Vert \mu_{f,\tau} \Vert =
\int_{\de \D} d\mu_{f,\tau} =\Re \frac{\tau+f(0)}{\tau-f(0)}.
\end{equation}
Moreover,   equality holds in (\ref{Cowen-Pomm-Ine}) if and
only if
$\sum_{j=1}^{n}\frac{1}{|f'(\zeta_j)|}\delta_{\zeta_j}=\mu_{f,\tau}$
if and only if
\begin{equation*}
\begin{split}
\Re \frac{\tau+f(z)}{\tau-f(z)}&=\int_{\de \D} P(\zeta,z)d\mu_{f,\tau}(\zeta)
= \sum_{j=1}^{n}\frac{1}{|f'(\zeta_j)|}\int_{\de \D} P(\zeta,z)d\delta_{\zeta_j}(\zeta) \\
&= \Re \sum_{j=1}^{n}\frac{1}{|f'(\zeta_j)|}\frac{\zeta_j+z}{\zeta_j-z}
\quad \hbox{for all\ } z\in \D,
\end{split}
\end{equation*}
namely, if and only if $f$ is a finite Blaschke product of
order $n$. Inequality (\ref{Cowen-Pomm-Ine}) was obtained in
\cite[Thm 8.1]{CoPo} by Cowen and Pommerenke with complete
different techniques.
\end{remark}

\section{Differentiability of Aleksandrov-Clark measures and the representation formula}

First of all we prove Proposition \ref{mainproposition}.

\begin{proof}[Proof of Proposition \ref{mainproposition}]
For the sake of simplicity, let us denote by
$\mu_t:=\mu_{t,\tau}$ the Aleksandrov-Clark measure of $\v_t$
at $\tau$. Moreover, for $t\geq 0$ we define
\[
\sigma_t:=\frac{\mu_t-\delta_\tau}{t}.
\]
Let $\sigma_t=\sigma_t^s+\sigma_t^a$ be the Lebesgue
decomposition of $\sigma_t$ with respect to the Lebesgue
measure. By Corollary \ref{univalegrupNev} it follows:
\begin{equation}\label{singular-part}
    \sigma_t^s=\frac{e^{-\lambda t}-1}{t}\delta_\tau.
\end{equation}
Taking the limit as $t\to 0$, we have
\[
\sigma_t^s\stackrel{w^\ast}{\longrightarrow} -\lambda
\delta_\tau.
\]
Now we examine the absolutely continuous part $\sigma_t^a$.
Since $\sigma_t^a\geq 0$ we have
\begin{equation*}
\begin{split}
\|\sigma_t^a\|&=\int_{\de\D} d\sigma_t^a=\int_{\de\D}
d\sigma_t-\int_{\de\D}d\sigma_t^s\stackrel{\eqref{singular-part}}{=}\int_{\de\D}
d\sigma_t-\frac{e^{-\lambda
t}-1}{t}\\&=\frac{1}{t}\int_{\de\D}d\mu_t-\frac{1}{t}-\frac{e^{-\lambda
t}-1}{t}
\stackrel{\eqref{AleCl}}{=}\frac{1}{t}\Re\left[\frac{\tau+\v_t(0)}{\tau-\v_t(0)}-1
\right]-\frac{e^{-\lambda t}-1}{t}\\
&=\Re\left[\frac{2\v_t(0)}{t}\frac{1}{\tau-\v_t(0)}
\right]-\frac{e^{-\lambda t}-1}{t},
\end{split}
\end{equation*}
and then, taking the limit for $t\to 0$ and by \eqref{genera},
we obtain
\begin{equation}\label{limt}
    \lim_{t\to 0}\|\sigma_t^a\|=2\Re G(0)+\lambda.
\end{equation}
This implies that $\{\|\sigma_t\|\}$ is uniformly bounded for
$t<<1$. Since the ball in the $\hbox{weak}^\ast$-topology of
measures on $\de\D$ is compact and metrizable, the net $\{\sigma_t\}$ is
sequentially compact.  Now, by \eqref{AleCl},
\begin{equation*}
\begin{split}
\int_{\de \D} P(\zeta,z) d\sigma_t(\zeta) &=\frac{1}{t}\Re
\left[
\frac{\tau+\v_t(z)}{\tau-\v_t(z)}-\frac{\tau+z}{\tau-z}\right]
\\&=2\Re \left[
\frac{\v_t(z)-z}{t}\cdot\frac{\tau}{(\tau-\v_t(z))(\tau-z)}\right],
\end{split}
\end{equation*}
and  \eqref{genera} yields
\begin{equation}\label{limstesso}
\lim_{t\to 0} \int_{\de \D} P(\zeta,z) d\sigma_t(\zeta)=   2\Re
\left[\frac{G(z)\tau}{(\tau-z)^2} \right].
\end{equation}
This implies that given two accumulation points $\sigma$ and
$\sigma'$ of $\{\sigma_t\}$ we have
$\int_{\de\D}P(\zeta,z)d\sigma(\zeta)=\int_{\de\D}P(\zeta,z)d\sigma'(\zeta)$.
Hence $\sigma=\sigma'$ (see, {\sl e.g.}, \cite[p. 10]{Koosis}).
Therefore the net $\{\sigma_t\}$ is actually
$\hbox{weak}^\ast$-convergent for $t\to 0$.

Finally, denote by $\mu$ the limit of $\{\sigma_t^a\}$. Since
$\sigma_t^a$ is a positive measure for all $t\geq 0$,   so is
$\mu$ and the proof  is completed.
\end{proof}

Now we are in the good shape to prove our representation
formula.

\begin{proof}[Proof of Theorem \ref{main}]
We retain the same notations as in the proof of Proposition
\ref{mainproposition}.

Since
\[
\Re \frac{\tau+\v_t(z)}{\tau-\v_t(z)}=\int_{\de \D} \Re\frac{\zeta+z}{\zeta-z}
d\mu_t(\zeta) \quad \hbox{for all\ } z\in \D,
\]
and $\frac{\tau+\v_t(z)}{\tau-\v_t(z)}$ and $\int_{\de \D}
\frac{\zeta+z}{\zeta-z} d\mu_t(\zeta)$ are analytic
functions, it follows
\[
\frac{\tau+\v_t(z)}{\tau-\v_t(z)}=\int_{\de \D} \frac{\zeta+z}{\zeta-z}
d\mu_t(\zeta) +i \Im \left( \frac{\tau+\v_t(0)}{\tau-\v_t(0)}\right) \quad \hbox{for all\ } z\in \D.
\]
Hence
\[
\frac{\tau+\v_t(z)}{\tau-\v_t(z)}-\frac{\tau+z}{\tau-z}=
\int_{\de \D} \frac{\zeta+z}{\zeta-z}d(\mu_t-\delta_\tau)(\zeta) +
i \Im \left( \frac{\tau+\v_t(0)}{\tau-\v_t(0)}-1\right) \quad \hbox{for all\ } z\in \D.
\]
After some computations and dividing by $t$ we obtain
\[
\frac{\v_t(z)-z}{t}=
\frac{(1-\v_t(z){\bar \tau})(\tau-z)}{2}
\left[
\int_{\de \D} \frac{\zeta+z}{\zeta-z}d\sigma_t(\zeta) +
2i \Im \left( \frac{\v_t(0)}{t}\frac{1}{\tau-\v_t(0)}\right)
\right]
\]
for all $z\in \D$. Now, by Proposition \ref{mainproposition},
passing to the limit as t goes to 0, we deduce
\[
G(z)=\frac{(1-z{\bar \tau})(\tau-z)}{2} \left[
-\lambda\frac{\tau+z}{\tau-z}+\int_{\de \D}
\frac{\zeta+z}{\zeta-z}d\mu(\zeta) + 2i \Im \left(
\frac{G(0)}{\tau}\right) \right].
\]
Setting  $p(z):=\frac{1}{2}\int_{\de \D}
\frac{\zeta+z}{\zeta-z}d\mu(\zeta) + i \Im \left(
\frac{G(0)}{\tau}\right)$, we obtain \eqref{descom-geninf}.

Moreover, since $\tau$ is a regular boundary fixed point, by
\cite[Theorem 1]{CoDiPo} it follows $\angle
\lim_{z\to\tau}\frac{G(z)}{z-\tau}=\lambda$. Then an easy
computation shows that $\angle \lim_{z\to \tau}
(z-\tau)p(z)=0.$

In order to prove uniqueness, assume that $q:\D\to\C$ is
holomorphic and $\gamma\in\C$ are such that $\Re q\geq 0$,
$\angle \lim _{z\to \tau}(z-\tau)q(z)=0$, and
\[
G(z)=({\bar
\tau}z-1)(z-\tau)\left[q(z)-\frac{\gamma}{2}\frac{\tau+z}{\tau-z}\right]
\quad\hbox{for all}\ z\in\D.
\]
Then,  by \cite[Theorem 1]{CoDiPo},
\[
\lambda =\angle \lim_{z\to\tau}\frac{G(z)}{z-\tau}=\angle
\lim_{z\to\tau} {\bar \tau}(z-\tau)
\left[q(z)-\frac{\gamma}{2}\frac{\tau+z}{\tau-z}\right]=\gamma
.
\]
From this, it follows immediately that $p=q$.

Now we prove the converse: let $p:\D\to\C$ holomorphic with
$\Re p\geq 0$, $\tau \in \de \D$ and $\lambda \in \R$, and let
 $G$ be defined by \eqref{descom-geninf}. We want to prove that $G$ is an infinitesimal
generator of some continuous semigroup of holomorphic self-maps
of the unit disc.

Since $\Re p \geq 0$, by \cite[Theorem 1.4.19]{Abate}, the
function
\[
H_1(z):=({\bar \tau}z-1)(z-\tau)p(z)  \quad\hbox{for all}\ z\in\D
\]
is the infinitesimal generator of a continuous semigroup of
holomorphic functions with Denjoy-Wolff point $\tau$. By
\cite[Theorem 1]{CoDiPo}, $\angle\lim_{z\to\tau}H_1(z)=0$ and
\[
\angle\lim_{z\to\tau}\frac{H_1(z)}{z-\tau}=\beta
\]
 for some $\beta\leq
0$. Our hypothesis that $\angle\lim_{z\to\tau} (z-\tau)p(z)=0$
implies that actually $\beta=0$. Therefore the semigroup
associated to $H_1$ has Denjoy-Wolff point $\tau$ with boundary
dilatation coefficient $1$ for all $t\geq 0$. In particular, if
$\lambda=0$ we are done.

Assume $\lambda\neq 0$. Then  $H_2:\D\to \C$ defined by
$H_2(z):=\lambda({\bar \tau}z-1)(z+\tau)$  is also the
infinitesimal generator of a semigroup of linear fractional
maps (in fact of hyperbolic automorphisms) with fixed points
$\tau$ and $-\tau$ (see, \cite[Corollary 1.4.16]{Abate}). Since
the set of infinitesimal generators is a real convex cone (see,
{\sl e.g.}, \cite[Corollary 1.4.15]{Abate}), it follows that
$G(z)=H_1(z)+H_2(z)$ is the infinitesimal generator of a
semigroup of holomorphic self-maps of the unit disc. Moreover,

\[
\angle \lim_{z\to\tau}G(z)=\angle \lim_{z\to\tau}H_1(z)+\angle \lim_{z\to\tau}H_2(z)=0
\]
and
\[
\angle \lim_{z\to\tau}\frac{G(z)}{z-\tau}=\angle
\lim_{z\to\tau}\frac{H_1(z)}{z-\tau}+\angle
\lim_{z\to\tau}\frac{H_2(z)}{z-\tau}=\beta+\lambda=\lambda.
\]
Therefore, by \cite[Theorem 1]{CoDiPo}, $\tau$ is a boundary
regular fixed point of the semigroup with boundary dilatation
coefficients $(e^{t\lambda})$.
\end{proof}

\section{Final Remarks}

\noindent{\bf 1.} The measure $\mu$ in formula
\eqref{maineq-intro} is strictly related to the infinitesimal
generator $G$ of $(\v_t)$. In fact, from classical measure
theory, if $\mu_r$ for $r\in (0,1)$ is the measure defined by
$\mu_r=\rho_r dm$  with density $\rho_r(\xi):=\int_{\de \D}
P(\zeta, r \xi)d\mu(\zeta)$, $\xi\in\de\D$, then
\begin{enumerate}
  \item $\lim_{r\to 1}\mu_r(\xi)=\mu^a(\xi)$ for $m$-almost
  every $\xi\in\de\D$ and
  \item $\mu_r\stackrel{w^\ast}{\longrightarrow} \mu$.
\end{enumerate}
From \eqref{maineq-intro} and \eqref{limstesso} it follows that
\[
\mu_r(\xi)=2\Re \left[ \frac{G(r\xi)\tau}{(\tau-r\xi)^2}
\right]dm(\xi)+\lambda \Re\left[ \frac{\tau+r\xi}{\tau-r\xi}
\right]dm(\xi).
\]
Thus, from (1) and (2) above we obtain
\begin{proposition}\label{Glimit}
Let $(\v_t)$ be a continuous semigroup of holomorphic self-maps
of the unit disc $\D$ with infinitesimal generator $G$. Let
$\tau\in\de\D$ be a boundary regular fixed point for $(\v_t)$.
Let $\mu$ be the positive measure defined in
\eqref{maineq-intro}. Then
\begin{itemize}
  \item[a)] $\Re \left[ \frac{G^\ast( \xi)\tau}{(\tau-
  \xi)^2}\right]\in L^1(\de \D, m)$ and $\mu^a(\xi)=2\Re \left[ \frac{G^\ast(\xi)\tau}{(\tau-\xi)^2}
\right]dm(\xi)$.
  \item[b)] $\int_{\de \D}f(\zeta)d\mu(\zeta)=
\lim_{r\to 1}\int_{\de \D}f(\zeta)\left(2\Re \left[ \frac{G(r\zeta)\tau}
  {(\tau-r\zeta)^2}\right]+\lambda \Re\left[ \frac{\tau+r\xi}{\tau-r\xi}\right]\right)dm(\zeta)$ for all $f\in C(\de\D)$.
\end{itemize}
\end{proposition}

\noindent{\bf 2.} From the proof of Theorem \ref{main} it
follows that the condition $\angle\lim_{z\to\tau}
(z-\tau)p(z)=0$ is not necessary in order to show that
\eqref{descom-geninf} defines an infinitesimal generator,
namely, what we really proved is:

\begin{proposition}
Let $p:\D\to\C$ holomorphic with $\Re p\geq 0$, $\tau \in \de
\D$ and $\lambda \in \R$. Then $\angle\lim_{z\to
\tau}(z-\tau)p(z)=\beta$ exists for some $\beta\leq 0$ and
 the function $G$ defined as in \eqref{descom-geninf}
is the infinitesimal generator of a semigroup of holomorphic
self-maps of the unit disc for which $\tau$ is a boundary
regular fixed point with boundary dilatation coefficients
$(e^{(\beta+\lambda) t})$.
\end{proposition}

\noindent{\bf 3.} Theorem \ref{main} shows that given a
semigroup of holomorphic functions $(\v_t)$ with a boundary
regular fixed point $\tau$, its infinitesimal generator $G$ is
the sum of the infinitesimal generator of a semigroup of
parabolic holomorphic maps with Denjoy-Wolff point at $\tau$
(namely, $H_1(z)=({\bar \tau}z-1)(z-\tau)p(z) $) plus, if
$\lambda \neq 0$, the infinitesimal generator of a group of
hyperbolic automorphisms of the unit disc (namely,
$H_2(z)=\frac{\lambda}{2}({\bar \tau}z^2-\tau)$) with a fixed
point at $\tau$. Notice that $\tau$ is the Denjoy-Wolff point
for  $(\v_t)$ if and only if it is the Denjoy-Wolff point for
the group of hyperbolic automorphisms.

\end{document}